\documentclass[reqno]{amsart}

\usepackage[margin=28mm]{geometry}
\usepackage{amsmath,amssymb,amsthm}
\usepackage[hidelinks]{hyperref}
\usepackage{color,cite}
\hypersetup{
  pdftitle={The uniform approximation property for Schatten classes},
  pdfsubject={The uniform approximation property in Schatten classes}
}
\newtheorem{theorem}{Theorem}
\newtheorem{corollary}[theorem]{Corollary}
\newtheorem{proposition}[theorem]{Proposition}
\allowdisplaybreaks[2]

\title{The uniform approximation property for Schatten classes}
\author[J. Huang]{Jinghao Huang}\thanks{Jinghao Huang was supported by the NNSF of China  (No.12031004, 12301160, 12471134 and 12671159).}
    \address{Institute for  Advanced Study in  Mathematics, Harbin Institute of Technology, Harbin, 150001, China}
\email{{\color{blue}jinghao.huang@hit.edu.cn}}   
    \author[R. Xu]{Ran Xu}
     \address{Institute for  Advanced Study in  Mathematics, Harbin Institute of Technology, Harbin, 150001, China}   \email{{\color{blue}xxurann@stu.hit.edu.cn}}

\begin{document}
\maketitle

\begin{abstract}
We show that the Schatten class \(S_p\), \(1\le p\le\infty\), has the
uniform approximation property if and only if \(p=2\). This resolves
a long-standing open problem in the field.


\end{abstract}

\section{Introduction}\label{sec:criterion}

Let \(S_p\) denote the Schatten class on \(\ell_2\), with
\(\|Z\|_p=(\operatorname{Tr}|Z|^p)^{1/p}\) for \(1\le p<\infty\).
We give \(S_\infty\), the space of compact operators on \(\ell_2\), the
operator norm and use \(1/\infty=0\).
A Banach space \(X\) has the \emph{uniform approximation property}
(UAP) \cite[p.~525]{jp} if there exist \(K\ge1\)
and \(r\colon\mathbb N\to\mathbb N\) such that, for every \(d\in\mathbb N\)
and every subspace \(E\subset X\) with \(\dim E\le d\), some linear \(T\colon X\to X\)
satisfies
\begin{equation}\label{eq:uap}
  T|_E=I_E,\qquad \|T\|\le K,\qquad \operatorname{rank}T\le r(d).
\end{equation}
Pe{\l}czy\'nski and Rosenthal proved that the classical spaces
\(L^p(\mu)\) have the UAP \cite{pr}. In 1976 Lindenstrauss asked
whether \(K(\ell_2)\) has the UAP
\cite[Question~8, first part, p.~18-05]{lind76}.
The case \(p=\infty\) of Theorem~\ref{thm:main} recovers Szankowski's
negative answer, as recorded by Pisier \cite[p.~176]{pisier95}.
In 1981 Szankowski proved that
\(B(\ell_2)\) fails even the ordinary approximation property (AP)
\cite[p.~89]{sz81}. This is the space denoted by \(S^\infty\) in
\cite[p.~1465]{px}. In 1984 he proved that \(S_p\) fails the UAP for
\(80<p<\infty\) \cite[p.~343]{sz84}. For \(1<p<\infty\), the two
conjugate ranges are \(p>80\) and \(p<80/79\)
\cite[Introduction]{lds}. In 2003 Pisier and Xu asked whether every
\(S_p\), \(p\ne2\), fails the UAP \cite[p.~1479]{px}. Within
\(1<p<\infty\), the exponents left by the 1984 bound are
\[
 [80/79,2)\ \cup\ (2,80].
\]

In 2004 Junge, Nielsen, Ruan and Xu \cite[Thm.~2.19]{jnrx} used Szankowski's
result to prove that for every \(80<p<\infty\) there is a finite von
Neumann algebra \(N\) with separable predual such that \(L^p(N)\)
fails the ordinary AP, that is, the identity cannot be approximated
uniformly on compact sets by finite-rank operators.
Corollary~\ref{cor:ap} extends their theorem to every
\(1<p<\infty\), \(p\ne2\). Corollary~\ref{cor:rigid-factor} realizes
any prescribed countable set of these exponents in a single
\(\mathrm{II}_1\) factor with property~\((T)\) and trivial outer
automorphism and fundamental groups. For each such \(p\),
Corollary~\ref{cor:olp} gives a Schatten ultrapower without the
ordinary bounded approximation property and applies this to
Problems~2.20(i) and (iv) of \cite[p.~285]{jnrx}.

For group von Neumann algebras, Lafforgue and de la Salle
\cite[Introduction and Thm.~A]{lds} constructed, for each \(p\ne2\),
discrete groups \(\Gamma\) whose group von Neumann algebras
\(L\Gamma\) yield spaces \(L^p(L\Gamma)\) without the completely
bounded approximation property (CBAP). De Laat and de la Salle
\cite[Thm.~1.1]{dlds} proved that one may take
\(\Gamma=\mathrm{SL}(n,\mathbb Z)\) for all sufficiently large \(n\)
depending on \(p\). The CBAP asks for finite-rank maps converging
pointwise in norm to the identity with uniformly bounded completely
bounded norms \cite[Sec.~2.2]{dlds}. On \(S_p\) itself, the orthogonal
projections \(P_n\) onto the first \(n\) coordinates give
\begin{equation}\label{eq:corner-cbap}
 C_n Z=P_nZP_n,\qquad \|C_n\|_{\mathrm{cb}}=1,\qquad
 \operatorname{rank}C_n=n^2,\qquad
 \lim_{n\to\infty}\|C_nZ-Z\|_p=0 \quad(Z\in S_p).
\end{equation}
Indeed, at matrix size \(k\) the amplification of \(C_n\) is compression
by \(P_n\otimes I_k\) on the amplified Schatten class and hence has
norm one \cite[Sec.~1.1, (2)]{lds}. Norm convergence in \eqref{eq:corner-cbap}
follows by density of the finite matrix corners, also in
\(S_\infty=K(\ell_2)\). Thus \(S_p\) has the CBAP for every \(p\).
Theorem~\ref{thm:main} concerns the exact restriction \(T|_E=I_E\)
and the bound \(r(\dim E)\) in \eqref{eq:uap}. For \(1\le p<\infty\),
\(p\ne2\), Corollary~\ref{cor:groups} below characterizes the countable
discrete groups for which \(L^p(L\Gamma)\) has this Banach-space property.

\begin{theorem}\label{thm:main}
For \(1\le p\le\infty\), the Schatten class \(S_p\) has the uniform
approximation property if and only if \(p=2\).
\end{theorem}

The proof occupies Section~\ref{sec:main-proof}. For \(p>2\) we use
the finite-matrix compressions and sign averaging
of Szankowski \cite[Secs.~2--4]{sz81}. The partitions in
Section~\ref{sec:fields} give the two singular-value lists in
\eqref{eq:local-spectra}. For \(2<p<\infty\), choosing the number
\(t\) of refinements with \(t(p-2)>8\) makes the bound in
\eqref{eq:level-bound} summable. The nuclear-norm estimate in
\eqref{eq:tail} is uniform in the terminal matrix size.
Corollary~\ref{cor:approximate-rank} quantifies the rank needed to
approximately fix a prescribed set of matrices. Theorem~\ref{thm:classification}
then classifies the Haagerup \(L^p\)-spaces with the UAP.
Section~\ref{sec:independent-row} treats the finite-dimensional
column-to-row norm independently of Theorem~\ref{thm:main}.

\section{Proof of Theorem~\ref*{thm:main}}\label{sec:main-proof}

\subsection{The Hilbert case and reduction to finite matrices}

For \(p=2\), the orthogonal projection onto \(E\subset S_2\) satisfies
\eqref{eq:uap} with \(K=1\) and \(r(d)=d\).

Write \(X_{p,D}=M_D(\mathbb C)\) with normalized Schatten norm
\(\|Z\|_{p,D}=D^{-1/p}\|Z\|_p\). For \(1/p+1/p'=1\), the trace
pairing of \cite[pp.~1464--1465]{px} gives
\begin{equation}\label{eq:pairing}
  \tau_D=D^{-1}\operatorname{Tr},\qquad
  \ell_A(Z)=\tau_D(A^*Z),\qquad
  \|\ell_A\|_{X_{p,D}^*}=\|A\|_{p',D}.
\end{equation}
If \eqref{eq:uap} holds on \(S_p\), it holds with the same \(K,r\)
on every \(X_{p,D}\). Indeed, the upper-left-corner embedding
\(Z\mapsto D^{-1/p}(Z\oplus0)\) is isometric and has the contractive
left inverse \(Y\mapsto D^{1/p}P_DYP_D\), where \(P_D\) is the
projection onto the first \(D\) coordinates. Composing \(T\) with
these two maps preserves the identity on the chosen subspace and
can only decrease its norm and rank.

Write \(\operatorname{tr}_X\) for the trace on endomorphisms of
\(X_{p,D}\), and \((\ell\otimes B)(Z)=\ell(Z)B\). For such an
endomorphism \(V\), let
\[
 \nu(V)=\inf_{V=\sum_i\ell_i\otimes B_i}
           \sum_i\|\ell_i\|_{X_{p,D}^*}\|B_i\|_{p,D}.
\]
Since \(\operatorname{tr}_X(T(\ell\otimes B))=\ell(TB)\), taking
the infimum gives
\begin{equation}\label{eq:trace-test}
  |\operatorname{tr}_X(TV)|\le\|T\|\nu(V).
\end{equation}
We use \eqref{eq:trace-test} in Section~\ref{sec:endpoints}.
Johnson and Pisier \cite[Prop.~1.2, (iii) implies (i)]{jp} use
the same trace decomposition for a general finite-rank sum.

\subsection{Finite-field partitions and Walsh submatrices}\label{sec:fields}

Fix \(2<p\le\infty\). Choose an integer \(t>4\) with \(t(p-2)>8\)
if \(p<\infty\). In \(\overline{\mathbb F}_2\), set
\(E_j=\mathbb F_{2^{2^{j+1}}}\) for \(0\le j\le t\), and fix
\(\lambda_j\in E_j\setminus E_{j-1}\) for \(1\le j\le t\).
Let \(L\) be a common multiple of \(t\) and \(2^{t+1}\). For
\(n\ge1\) and an integer \(s\ge0\), put
\begin{equation}\label{eq:parameters}
 m_n=2^{L(n+s)},\qquad F_n=\mathbb F_{m_n},\qquad
 q_n=m_n^2,\qquad k_n=q_{n+1}^2,\qquad b_n=k_n^{1/t}=m_{n+1}^{4/t}.
\end{equation}
Since \(2^{t+1}\mid L(n+s)\), the subfield criterion
\cite[Sec.~2.1]{ln} gives \(E_t\subseteq F_n\).
Since \(t\mid L\), \(b_n\) is an integer. By
\eqref{eq:parameters}, \(b_n/m_n=2^{4L/t}m_n^{4/t-1}\).
Choose \(s\) so that \(b_1/m_1\le1/2\). Since \(t>4\), this gives
\(b_n\le m_n/2\) for every \(n\).

Fix \(n\), write \(F=F_n\), \(m=m_n\), \(q=q_n\), \(b=b_n\),
and put \(V=F^2\). Index the \(k_n=b^t\) words by
\(h=(d_1,\ldots,d_t)\in\{1,\ldots,b\}^t\).
At level \(j\), a parent \(P\) consists of words with a fixed prefix
of length \(j-1\). Its children \(G\) are obtained by fixing \(d_j\).
Let \(\mathcal L_j\) be the partition of \(V\) into the lines
\(\{(x,\lambda_jx+w):x\in F\}\), \(w\in F\).

For \(Q=|E_{j-1}|\), let \(A_a(x)=x+ax^Q\), \(a\in F\).
These maps are \(E_{j-1}\)-linear. Since \(Q-1\mid m-1\) and
\(F^\times\) is cyclic \cite[Sec.~2.1]{ln},
\[
 \{a:\ker A_a\ne\{0\}\}=\{x^{1-Q}:x\in F^\times\},
 \qquad \#\{a:\ker A_a=\{0\}\}
 =m-\frac{m-1}{Q-1}\ge\frac m2\ge b .
\]
Here \(Q\ge4\). Assign a distinct invertible \(A_a\) to each digit
\(d\) at position \(j\), put \(D_{j,d}(x,y)=(A_a(x),A_a(y))\), and set
\begin{equation}\label{eq:permutation}
 \pi_h=D_{t,d_t}\circ\cdots\circ D_{1,d_1}.
\end{equation}
For \(\ell>j\), the coordinate map defining \(D_{\ell,d}\) is
\(E_j\)-linear, so \(A_a(y)-\lambda_jA_a(x)=A_a(y-\lambda_jx)\).
Thus \(D_{\ell,d}\) permutes the lines in \(\mathcal L_j\). By
\eqref{eq:permutation}, every \(h\in G\) induces the same partition
\[
 \mathcal A_G=\pi_h^{-1}\mathcal L_j
 =(D_{j,d_j}\circ\cdots\circ D_{1,d_1})^{-1}\mathcal L_j .
\]

Let \(G,G'\) be distinct children of \(P\), with digit parameters
\(a,a'\). After applying the common prefix permutation, intersections
of blocks in \(\mathcal A_G,\mathcal A_{G'}\) are fibers of the additive map
\[
 (x,y)\longmapsto
 \bigl(A_a(y)-\lambda_j A_a(x),\, A_{a'}(y)-\lambda_j A_{a'}(x)\bigr).
\]
In its kernel, subtraction gives
\((a-a')(y^Q-\lambda_jx^Q)=0\), hence \(y^Q=\lambda_jx^Q\) and
\(y=\lambda_jx\). Since \(\lambda_j\notin\mathbb F_Q\),
\((\lambda_j^Q-\lambda_j)x^Q=0\) implies \(x=y=0\).
The map is therefore a bijection of \(F^2\), and
\begin{equation}\label{eq:incidence}
 |A\cap A'|=1
 \qquad(A\in\mathcal A_G,\ A'\in\mathcal A_{G'},\ G\ne G'
        \text{ children of }P).
\end{equation}

With \(\chi(u)=(-1)^{\operatorname{Tr}_{F/\mathbb F_2}(u)}\), define
the real \(q\times q\) matrices
\begin{equation}\label{eq:walsh}
 V_h((x,y),z)=\chi(xv+yu),\qquad (u,v)=\pi_h(z),\qquad
 V_h^*V_h=qI_q .
\end{equation}
The character orthogonality identity
\(\sum_{x\in F}\chi(cx)=m\mathbf 1_{\{c=0\}}\)
\cite[Secs.~2.3 and 5.1]{ln}, applied in both coordinates, gives
the last equality in \eqref{eq:walsh} since \(\pi_h\) is bijective.
Szankowski's submatrix estimates \cite[Lemma~1, (2.1)--(2.2)]{sz81}
distinguish two cases. We compute the singular values at every digit \(j\).
Fix \(S=\{(x,\lambda_jx+w):x\in F\}\in\mathcal L_j\) and
\(A\in\mathcal A_G\), and let \(V_h[S,A]\) be the submatrix with
row set \(S\) and column set \(A\). For \(z,z'\in A\), write
\((u,v)=\pi_h(z)-\pi_h(z')\). By \eqref{eq:walsh},
\begin{equation}\label{eq:walsh-gram}
 \bigl(V_h[S,A]^*V_h[S,A]\bigr)_{z,z'}
 =\chi(wu)\sum_{x\in F}\chi\bigl(x(v+\lambda_j u)\bigr)
 =m\chi(wu)\mathbf 1_{\{v=\lambda_j u\}} .
\end{equation}
If \(h\in G\), then \(\pi_h(A)\in\mathcal L_j\).
Writing \(\pi_h(z)=(u_z,v_z)\), \eqref{eq:walsh-gram} gives the
rank-one Gram matrix  $$ m(\chi(wu_z)\chi(wu_{z'}))_{z,z'\in A} $$
of trace \(m^2\). If \(h\in P\setminus G\), let \(G'\) be the
child containing \(h\). By \eqref{eq:incidence},
\(|A\cap\pi_h^{-1}(S')|=1\) for every \(S'\in\mathcal L_j\).
Thus \(v=\lambda_j u\) in \eqref{eq:walsh-gram} holds exactly when
\(z=z'\), and the Gram matrix is \(mI_m\). Writing \(s_i\) for
the singular values in decreasing order, we obtain
\begin{equation}\label{eq:local-spectra}
 \bigl(s_1(V_h[S,A]),\ldots,s_m(V_h[S,A])\bigr)=
 \begin{cases}
 (m,0,\ldots,0),& h\in G,\\
 (\sqrt m,\ldots,\sqrt m),& h\in P\setminus G .
 \end{cases}
\end{equation}

\subsection{Finite orthogonal matrices and their rectangular compressions}
\label{sec:matrices}

Put \(J_n=F_n^2\), \(K_n=J_n\times J_n\), \(D_n=\prod_{\ell=1}^nK_\ell\),
and \(N_n=|D_n|=\prod_{\ell=1}^n q_\ell^2\). Take \(D_0\) to be a
singleton. By \eqref{eq:parameters}, \(|K_{n+1}|=k_n=b_n^t\).
Identify \(K_{n+1}\) with the length-\(t\) label words of
Section~\ref{sec:fields}, and write \(V_h^{(n)}\) for its matrices at
level \(n\), with rows and columns indexed by \(J_n\).

Fix \(M\ge2\), put \(D=N_M\), and set \(U_{M+1,M}=(1)\).
For \(1\le n\le M\), rows and columns of \(U_{n,M}\) are indexed by
\(((x,y),i)\), where \(x,y\in J_n\) and
\(i\in\prod_{\ell=n+1}^M K_\ell\). Write \(i_1\) for the first
entry of \(i\). At \(n=M\), take \(i_1=*\) for the empty tail and
let \(V_*^{(M)}\) be any one of the matrices in \eqref{eq:walsh} at
level \(M\). Following the two-factor formula of
\cite[p.~96]{sz81}, define
\begin{equation}
 U_{n,M}((x,y),i;(u,v),j)
  =q_n^{-1}V_{i_1}^{(n)}(v,x)V_{j_1}^{(n)}(y,u)U_{n+1,M}(i,j).
 \label{eq:gate}
\end{equation}
Assume \(U_{n+1,M}\) is symmetric and orthogonal. Interchanging
\(((x,y),i)\) and \(((u,v),j)\) in \eqref{eq:gate} gives symmetry.
Since \(V_{j_1}^{(n)}(V_{j_1}^{(n)})^{\mathsf T}=q_nI\) by
\eqref{eq:walsh}, summing over \(u\) first yields
\begin{align*}
 &(U_{n,M}U_{n,M}^{\mathsf T})_{((x,y),i),((x',y'),i')}\\
 &\quad=\frac{\delta_{yy'}}{q_n}
  \left(\sum_j U_{n+1,M}(i,j)U_{n+1,M}(i',j)\right)
  \sum_{v\in J_n} V_{i_1}^{(n)}(v,x)V_{i'_1}^{(n)}(v,x')\\
 &\quad=\delta_{yy'}\delta_{ii'}\delta_{xx'}.
\end{align*}
Backward induction and \eqref{eq:gate} therefore give
\begin{equation}\label{eq:flat}
 U=U_{1,M}=U^{\mathsf T},\qquad U^{\mathsf T}U=I,\qquad
 |U(i,j)|=\prod_{\ell=1}^M q_\ell^{-1}=D^{-1/2}\quad(i,j\in D_M).
\end{equation}

Fix \(1\le n<M\), and write \(N=N_n\), \(m=m_n\),
\(k=k_n\), \(b=b_n\). Fix \(\alpha,\beta\in D_{n-1}\),
\(e,f\in J_n\), a child \(G\) of \(P\) at refinement \(j\),
\(S\in\mathcal L_j\), and \(A\in\mathcal A_G\). For
\(H\subseteq K_{n+1}\), let
\[
 \begin{aligned}
 E&=\{(\alpha,(e,s),i):s\in S,\ i\in\textstyle\prod_{\ell=n+1}^M K_\ell\},\\
 F_H&=\{(\beta,(u,f),j'):u\in A,\ j'\in\textstyle\prod_{\ell=n+1}^M K_\ell,\
                       j'_1\in H\}.
 \end{aligned}
\]
Set \(Z_H=P_EUP_{F_H}\), where the two projections select the indicated
coordinates. Let \(P_h\) be the projection on the tail coordinates
with first entry \(h\). Write \(\sim\) for equality after row and column
permutations, row and column sign changes, and insertion or deletion
of zero rows and columns. With \(V_h[S,A]\) taken at level \(n\),
\begin{equation}
 Z_H\ \sim\ N^{-1/2}\sum_{h\in H}V_h[S,A]\otimes(U_{n+1,M}P_h).
 \label{eq:rectangle}
\end{equation}
Indeed, in \eqref{eq:gate} all signs at levels \(\ell\le n-2\)
are fixed by \(\alpha,\beta\). For \(n\ge2\), write
\(\alpha_{n-1}=(x,y)\), \(\beta_{n-1}=(v,w)\). The two signs at
level \(n-1\) and the first sign at level \(n\) have product
\[
 \underbrace{V_{(e,s)}^{(n-1)}(w,x)V_{i_1}^{(n)}(f,e)}_{\text{row sign}}
 \ \underbrace{V_{(u,f)}^{(n-1)}(y,v)}_{\text{column sign}} .
\]
At \(n=1\) the only such factor is \(V_{i_1}^{(1)}(f,e)\).
After these sign changes, the entry indexed by \((s,i),(u,j')\) is
\[
 \left(\prod_{\ell=1}^n q_\ell^{-1}\right)
 V_{j'_1}^{(n)}(s,u)\, U_{n+1,M}(i,j')\,\mathbf 1_H(j'_1).
\]
The product equals \(N^{-1/2}\), proving \eqref{eq:rectangle}.

To compute the singular values in \eqref{eq:rectangle}, left multiply
by the orthogonal matrix \(I\otimes U_{n+1,M}^{\mathsf T}\). This gives
\begin{equation}\label{eq:direct-sum}
 \begin{aligned}
 (I\otimes U_{n+1,M}^{\mathsf T})
 \sum_{h\in H}V_h[S,A]\otimes(U_{n+1,M}P_h)
 &=\sum_{h\in H}V_h[S,A]\otimes P_h\\
 &\sim \bigoplus_{h\in H}
          \bigl(V_h[S,A]\otimes I_{D/(Nk)}\bigr).
 \end{aligned}
\end{equation}
The last step groups the tail coordinates by their first entry and
deletes the zero blocks. Each \(h\in K_{n+1}\) occurs in
\(\prod_{\ell=n+2}^M|K_\ell|=D/(Nk)\) tail coordinates. Thus
each singular value of \(N^{-1/2}V_h[S,A]\) contributes \(D/(Nk)\)
copies to \eqref{eq:rectangle}, counting nonzero values with multiplicity.
Writing \(c=|G|\), we have \(|P\setminus G|=(b-1)c\).
By \eqref{eq:local-spectra}, \eqref{eq:direct-sum} and the
normalized trace on \(D\times D\) matrices,
\begin{align}
 \|Z_G\|_{p',D}
   &=N^{-1/2}m\left(\frac{c}{Nk}\right)^{1/p'},\nonumber\\
 \|Z_{P\setminus G}\|_{p,D}
   &=N^{-1/2}\sqrt m\left(\frac{m(b-1)c}{Nk}\right)^{1/p}.
 \label{eq:rect-norms}
\end{align}
For \(p=\infty\), the second identity uses \(a^{1/\infty}=1\) for \(a>0\).

\subsection{Nuclear norms of the successive differences}\label{sec:increments}

For partitions \(\mathcal R,\mathcal C\) of \(D_M\), write \(P_A\)
for the coordinate projection onto \(A\subset D_M\), and set
\begin{equation}\label{eq:blocks}
 X_{A,B}=P_AUP_B,\qquad
 \mathcal B(\mathcal R,\mathcal C)
  =\sum_{A\in\mathcal R,\,B\in\mathcal C}\ell_{X_{A,B}}\otimes X_{A,B}.
\end{equation}
Let \(\mathcal P_n\) partition \(D_M\) by the first \(n\) coordinates,
and put \(\beta_n=\mathcal B(\mathcal P_n,\mathcal P_n)\).
The rectangles \(A\times B\) partition \(D_M\times D_M\). Thus
\eqref{eq:pairing} and \eqref{eq:blocks} give
\begin{equation}\label{eq:beta}
 \operatorname{tr}_X\beta_n
  =\sum_{A,B}\tau_D(X_{A,B}^*X_{A,B})=\tau_D(U^*U)=1,\qquad
 \operatorname{rank}\beta_n\le N_n^2,\qquad \beta_M=D^{-2}I.
\end{equation}
Here \(I\) acts on \(X_{p,D}\). The rank bound counts the \(N_n^2\)
terms of \eqref{eq:blocks}. For the last identity in \eqref{eq:beta},
each singleton rectangle gives \(X_{\{i\},\{j\}}=\pm D^{-1/2}e_{ij}\) by \eqref{eq:flat},
so its contribution at \(Z\) equals \(D^{-2}Z_{ij}e_{ij}\).

As in the row and column decomposition of \cite[(3.0)--(3.1)]{sz81},
we change one partition at a time. To pass from \(\beta_n\) to
\(\beta_{n+1}\), first refine the
column partition through the \(t\) digits of \(K_{n+1}\), keeping
the row partition at \(\mathcal P_n\). Then refine the row partition
through the same \(t\) digits, keeping the columns at \(\mathcal P_{n+1}\).
At one digit, \(P\) ranges over the parent sets of labels, and
\(G\subset P\) in a sum ranges over their immediate children.

\paragraph{Column refinement.}
For \(a,d\in D_n\), let \(X_{a;d,H}\) be \(U\) restricted to rows
with prefix \(a\) and columns with prefix \(d\) and next label in
\(H\), extended by zero. Substituting
\(X_{a;d,P}=\sum_{G\subset P}X_{a;d,G}\) in \eqref{eq:blocks},
the operator before the refinement minus the operator after it is
\begin{equation}\label{eq:split}
 \sum_P\sum_{G\subset P}\sum_{a,d\in D_n}
       \ell_{X_{a;d,G}}\otimes X_{a;d,P\setminus G}.
\end{equation}
For fixed \(G\), partition the prefix pairs in \eqref{eq:split}
by the rectangles of \eqref{eq:rectangle}. Within one rectangle,
\(a=(\alpha,(e,s))\), \(d=(\beta,(u,f))\), and
\((s,u)\in S\times A\). Take independent uniform signs
\((\varepsilon_a)_{a\in D_n},(\eta_d)_{d\in D_n}\).
The sign-averaging identity of \cite[Sec.~3, (3.6)--(3.9)]{sz81}
applies since
\(\mathbb E[\varepsilon_a\eta_d\varepsilon_{a'}\eta_{d'}]
=\mathbf 1_{a=a'}\mathbf 1_{d=d'}\). With all sums below over this
rectangle, it gives
\begin{equation}\label{eq:sign-average}
 \sum_{a,d}\ell_{X_{a;d,G}}\otimes X_{a;d,P\setminus G}
 =\mathbb E\!\left[
    \ell_{\sum_{a,d}\varepsilon_a\eta_d X_{a;d,G}}
      \otimes\sum_{a,d}\varepsilon_a\eta_d X_{a;d,P\setminus G}
               \right].
\end{equation}
For \(H=G,P\setminus G\), the signed sum is \(Z_H\) multiplied
on the left and right by diagonal unitaries, with entries
\(\varepsilon_a\) on row prefix \(a\) and \(\eta_d\) on column
prefix \(d\). By \eqref{eq:pairing} and \eqref{eq:sign-average}, the
nuclear norm of this rectangle's contribution is at most
\(\|Z_G\|_{p',D}\|Z_{P\setminus G}\|_{p,D}\).

\paragraph{Row refinement.}
Let \(X_{a,H;d,h}\) be \(U\) restricted to rows with prefix \(a\)
and next label in \(H\), and columns with prefix \((d,h)\), extended
by zero. Expanding \eqref{eq:blocks} as in \eqref{eq:split} gives
the terms \(\ell_{X_{a,G;d,h}}\otimes X_{a,P\setminus G;d,h}\).
In a transposed rectangle take \(a=(\beta,(u,f))\),
\(d=(\alpha,(e,s))\), \((s,u)\in S\times A\), and all \(h\in K_{n+1}\).
For independent uniform signs indexed by \(a\) and \((d,h)\),
\[
 \mathbb E[\varepsilon_a\eta_{d,h}\varepsilon_{a'}\eta_{d',h'}]
 =\mathbf 1_{a=a'}\mathbf 1_{d=d'}\mathbf 1_{h=h'},\qquad
 \sum_{a,d,h} X_{a,H;d,h}=Z_H^{\mathsf T}.
\]
The second identity follows from the symmetry in \eqref{eq:flat} and the union
over all \(h\) of the columns with prefix \((d,h)\). Applying
\eqref{eq:sign-average} with \(d\) replaced by \((d,h)\) bounds the
nuclear norm of this contribution by
\(\|Z_G^{\mathsf T}\|_{p',D}\|Z_{P\setminus G}^{\mathsf T}\|_{p,D}
=\|Z_G\|_{p',D}\|Z_{P\setminus G}\|_{p,D}\).

At one digit all children have size \(c\). For each of the \(k/c\)
children, the \(N^2\) prefix pairs split into \(N^2/m^2\) rectangles,
each containing \(m^2\) pairs. By \eqref{eq:rect-norms},
\begin{equation}\label{eq:one-step}
 \begin{aligned}
 \|Z_G\|_{p',D}\|Z_{P\setminus G}\|_{p,D}
  &=\frac{c}{kN^2}\,m^{3/2+1/p}(b-1)^{1/p},\\
 \frac{k}{c}\frac{N^2}{m^2}\,
 \|Z_G\|_{p',D}\|Z_{P\setminus G}\|_{p,D}
  &=m^{-1/2+1/p}(b-1)^{1/p}.
 \end{aligned}
\end{equation}
Summing \eqref{eq:one-step} over the \(t\) column refinements and
the \(t\) row refinements yields
\begin{equation}\label{eq:level-bound}
 \nu(\beta_n-\beta_{n+1})
 \le 2t\,m_n^{-1/2+1/p}(b_n-1)^{1/p}\qquad(1\le n<M).
\end{equation}

\subsection{The trace contradiction and the conjugate exponents}\label{sec:endpoints}

For \(2<p\le\infty\), set \(\delta=1/2-1/p-4/(tp)>0\), with
\(\delta=1/2\) at \(p=\infty\). From \eqref{eq:parameters},
\(b_n=(2^L m_n)^{4/t}\). By \eqref{eq:beta} and \eqref{eq:level-bound},
for every \(M>n_0\) we have
\begin{equation}\label{eq:tail}
 \begin{aligned}
 \nu(\beta_{n_0}-N_M^{-2}I)
 &\le\sum_{n=n_0}^{M-1}\nu(\beta_n-\beta_{n+1})\\
 &\le 2t\,2^{4L/(tp)}\sum_{n=n_0}^\infty m_n^{-\delta}
  =\frac{2t\,2^{4L/(tp)}}{1-2^{-L\delta}}\,m_{n_0}^{-\delta}.
 \end{aligned}
\end{equation}
Suppose \(S_p\) has the UAP with \(K,r\) as in \eqref{eq:uap}.
Choose \(n_0\) so that the right-hand side of \eqref{eq:tail} is
less than \(1/(2K)\). Put \(d=N_{n_0}^2\), then choose \(M>n_0\)
with \(N_M^2>2Kr(d)\). By \eqref{eq:beta}, the range of
\(\beta_{n_0}\) has dimension at most \(d\). Apply \eqref{eq:uap}
in \(X_{p,N_M}\) to obtain \(T\) with \(T\beta_{n_0}=\beta_{n_0}\).
The sum of the eigenvalues of \(T\), counted with algebraic
multiplicity, satisfies
\(|\operatorname{tr}_X T|\le\operatorname{rank}(T)\|T\|\).
It follows from \eqref{eq:beta}, \eqref{eq:trace-test} and
\eqref{eq:tail} that
\begin{equation}\label{eq:contradiction}
 \begin{aligned}
 1=|\operatorname{tr}_X(T\beta_{n_0})|
 &\le \|T\|\,\nu(\beta_{n_0}-N_M^{-2}I)
      + N_M^{-2}|\operatorname{tr}_X T|\\
 &< \tfrac12+N_M^{-2}Kr(d)<1 .
 \end{aligned}
\end{equation}

For \(1\le u<2\), take \(p=u'\), including \(p=\infty\) when \(u=1\).
Under the conjugate-linear isometry \(A\mapsto\ell_A\) in
\eqref{eq:pairing}, the Banach adjoint on \(X_{u,D}\) satisfies
\[
 (\ell_A\otimes B)^*=\ell_B\otimes A,\qquad
 \|\ell_B\|_{X_{u,D}^*}\|A\|_{u,D}=\|B\|_{p,D}\|A\|_{p',D}.
\]
Taking infima over decompositions and applying the same identity
twice shows that adjunction preserves the nuclear norm. Each term
of \(\beta_n\) in \eqref{eq:blocks} has \(A=B\) in the adjoint
formula, so \(\beta_n^*=\beta_n\) under this identification.
Consequently \eqref{eq:beta} and \eqref{eq:tail} hold on \(X_{u,D}\)
for the same matrices and the parameters chosen for \(p=u'\).
If \(S_u\) had the UAP, \eqref{eq:uap} on \(X_{u,N_M}\) would
again give \eqref{eq:contradiction}, with \(K,r\) chosen for \(S_u\).
In particular, for \(1\le p\le\infty\), \(p\ne2\), no common \(K,r\)
in \eqref{eq:uap} work for the family \((X_{p,D})_{D\ge1}\), since
such a pair would give \eqref{eq:contradiction} on \(X_{p,N_M}\).
This completes the proof of Theorem~\ref{thm:main}. 

\section{Applications}\label{sec:applications}

\subsection{The rank required to approximately fix finitely many matrices}

Here \(S_p^D=M_D(\mathbb C)\) carries the usual Schatten norm.

\begin{corollary}\label{cor:approximate-rank}
Fix \(1\le p\le\infty\), \(p\ne2\), and \(\varepsilon>0\). There is an
integer \(m=m(p,\varepsilon)\) such that, for arbitrarily large \(D\),
there are \(x_1,\ldots,x_m\in S_p^D\) with \(\|x_j\|_p=1\) and the
following property. For every \(K>0\), every \(\eta\ge0\), and every
linear \(T:S_p^D\to S_p^D\),
\begin{equation}\label{eq:approximate-rank}
 \|T\|\le K,\qquad \max_{1\le j\le m}\|Tx_j-x_j\|_p\le\eta
 \quad\Longrightarrow\quad
 \frac{\operatorname{rank}T}{D^2}\ge\frac{1-\eta}{K}-\varepsilon .
\end{equation}
\end{corollary}

\begin{proof}
Use the parameters of Section~\ref{sec:fields} for \(\max\{p,p'\}\).
By \eqref{eq:tail} and the adjoint identity for \(\beta_n\) in
Section~\ref{sec:endpoints}, choose \(n\) such that, on \(X_{p,D}\),
\[
 \nu(\beta_n-D^{-2}I)\le\varepsilon
 \qquad(M>n,\ D=N_M).
\]
Put \(m=N_n^2\). For every \(A,B\in\mathcal P_n\), use the test
matrix \(x_{A,B}=X_{A,B}/\|X_{A,B}\|_p\) from \eqref{eq:blocks}.

To compute its norm, fix the two prefixes defining \(A,B\) and
expand \eqref{eq:gate} through level \(n\). At levels \(1,\ldots,n-1\)
both signs are fixed. At level \(n\) the first sign depends only on
the row tail and the second only on the column tail. Thus the
\(A\times B\) block is \(N_n^{-1/2}U_{n+1,M}\) up to row and column
signs. Since \(|B|=D/N_n\) and \(U_{n+1,M}\) is orthogonal,
\begin{equation}\label{eq:prefix-singular}
 X_{A,B}^*X_{A,B}=N_n^{-1}P_B,\qquad
 \|X_{A,B}\|_{q,D}=N_n^{-1/2-1/q}\quad(1\le q\le\infty).
\end{equation}
The \(N_n^2\) pairs in \eqref{eq:prefix-singular} therefore satisfy
\[
 \sum_{A,B}\|X_{A,B}\|_{p',D}\|X_{A,B}\|_{p,D}
 =\sum_{A,B}\tau_D(X_{A,B}^*X_{A,B})=N_n^2 N_n^{-2}=1 .
\]
Since \(\|Z\|_{p,D}=D^{-1/p}\|Z\|_p\), the operator norm of \(T\)
is the same on \(S_p^D\) and \(X_{p,D}\). The test inequalities give
\(\|(I-T)X_{A,B}\|_{p,D}\le\eta\|X_{A,B}\|_{p,D}\).
By \eqref{eq:pairing} and \eqref{eq:blocks},
\[
 \bigl|\operatorname{tr}_X((I-T)\beta_n)\bigr|
 \le \sum_{A,B}\|X_{A,B}\|_{p',D}\|(I-T)X_{A,B}\|_{p,D}
 \le\eta .
\]
Since \(\operatorname{tr}_X\beta_n=1\) by \eqref{eq:beta},
\eqref{eq:trace-test} yields
\[
 1-\eta\le|\operatorname{tr}_X(T\beta_n)|
 \le \|T\|\nu(\beta_n-D^{-2}I)+D^{-2}|\operatorname{tr}_X T|
 \le K\varepsilon+K\,\frac{\operatorname{rank}T}{D^2}.
\]
Division by \(K\) gives \eqref{eq:approximate-rank}.
\end{proof}

For \(K=1\) and \(\eta=\varepsilon=0.01\), \eqref{eq:approximate-rank}
requires \(\operatorname{rank}T\ge0.98D^2\). If \(T=BA\) with linear
maps \(A:S_p^D\to\mathbb C^r\) and \(B:\mathbb C^r\to S_p^D\), then
\(r\ge\operatorname{rank}T\). Thus the same lower bound applies to
the number of complex linear measurements.

\subsection{Classification of Haagerup \(L^p\)-spaces}\label{sec:classification}

Let \(\mathcal M\) be a von Neumann algebra and let \(L^p(\mathcal M)\)
denote its Haagerup \(L^p\)-space.
The algebra \(\mathcal M\) is \emph{subhomogeneous} if there is a
finite integer \(s\) such that
\begin{equation}\label{eq:subhomogeneous}
 \mathcal M\simeq\bigoplus_{k=1}^{s}
      L^\infty(\Omega_k,\mu_k)\,\overline{\otimes}\,M_k .
\end{equation}
Zero summands are allowed \cite[Sec.~2]{lmz}.

\begin{theorem}\label{thm:classification}
For every von Neumann algebra \(\mathcal M\) and every
\(1\le p<\infty\), \(p\ne2\), the Banach space \(L^p(\mathcal M)\)
has the UAP if and only if \(\mathcal M\) is subhomogeneous.
\end{theorem}

\begin{proof}
First suppose \(\mathcal M\) is semifinite and choose a faithful normal
semifinite trace \(\tau\). The Haagerup space is isometric to the
tracial space \(L^p(\mathcal M,\tau)\) by \cite[Thm.~3.4(i)]{px}.
If \(\mathcal M\) is subhomogeneous, choose the measures in
\eqref{eq:subhomogeneous} so that \(\tau\) integrates the ordinary
matrix traces. Let \(S_p^k\) denote \(M_k\) with
the Schatten \(p\)-norm. Taking matrix entries gives
\[
 L^p(\mathcal M,\tau)
 =\left(\bigoplus_{k=1}^{s}L^p(\Omega_k;S_p^k)\right)_{\ell_p}
 \simeq
 \left(\bigoplus_{k=1}^{s} L^p(\Omega_k;\ell_p^{\,k^2})\right)_{\ell_p}.
\]
The entry maps and their inverses have uniformly bounded norms since
\(k\le s\). The space on the right is scalar \(L^p\) on the
disjoint union of \(k^2\) copies of each \(\Omega_k\). It has the UAP
by \cite{pr}, and conjugating the operators in \eqref{eq:uap} by the
displayed isomorphism transfers the UAP to \(L^p(\mathcal M,\tau)\).

Suppose \(\mathcal M\) is not subhomogeneous. For each \(D\ge2\),
choose matrix units \((v_{ij})_{i,j=1}^D\) in \(\mathcal M\)
\cite[proof of Lemma~2.1]{lmz}. By semifiniteness there is a projection
\(0\ne f\le v_{11}\) with \(\tau(f)<\infty\). The operators
\(v_{i1}fv_{1j}\) are matrix units and their diagonal sum \(e\) has
trace \(D\tau(f)\). Let \(\iota:M_D\to e\mathcal M e\) be the resulting
unital embedding. The trace-preserving conditional expectation
\(\mathbb E:e\mathcal M e\to\iota(M_D)\) acts contractively on
\(L^p\) \cite[Sec.~7, pp.~1490--1491]{px}. Since
\(\tau\circ\iota=\tau(e)\tau_D\), the maps
\begin{equation}\label{eq:matrix-retraction}
 J(a)=\tau(e)^{-1/p}\iota(a),\qquad
 Q(x)=\tau(e)^{1/p}\iota^{-1}\bigl(\mathbb E(exe)\bigr)
\end{equation}
satisfy \(QJ=I_{X_{p,D}}\) and \(\|J\|=\|Q\|=1\). If \(L^p(\mathcal M,\tau)\)
had UAP data \(K,r\), then for each \(E\subset X_{p,D}\) we could apply
\eqref{eq:uap} to \(J(E)\). The operator \(QTJ\) obtained from
\eqref{eq:matrix-retraction} would fix \(E\), have norm at most \(K\)
and rank at most \(r(\dim E)\). This would give the same \(K,r\) for
all \(D\). Apply Corollary~\ref{cor:approximate-rank} with
\(\eta=0\), \(\varepsilon=1/(2K)\) and \(E\) the span of its \(m\)
test matrices. Then \(D^2/(2K)\le\operatorname{rank}(QTJ)\le r(m)\)
for arbitrarily large \(D\), a contradiction. This proves the
classification for every semifinite algebra.

Every type~I algebra is semifinite. If \(\mathcal M\) is not of
type~I, let \(R\) be the hyperfinite \(\mathrm{II}_1\) factor.
By the semifinite case, \(L^p(R)\) fails the UAP. By
\cite[Thm.~3.5, p.~1476]{px}, \(L^p(\mathcal M)\) contains an
isometric, \(1\)-complemented copy of \(L^p(R)\) for the full range
\(1\le p<\infty\). Thus there are
contractions \(J:L^p(R)\to L^p(\mathcal M)\) and
\(Q:L^p(\mathcal M)\to L^p(R)\) with \(QJ=I\). If \(L^p(\mathcal M)\)
had the UAP, applying \eqref{eq:uap} to \(J(E)\) and replacing
\(T\) by \(QTJ\) would give the UAP on \(L^p(R)\). This contradiction
proves the theorem.
\end{proof}

In particular \(L^p(\mathcal M)\) fails the UAP for every type~III
factor and every \(1\le p<\infty\), \(p\ne2\). For semifinite
\((\mathcal M,\tau)\) in this range, combining Theorem~\ref{thm:classification}
with Ray and Sarkar's classification \cite[Thm.~1.1]{rs} gives
\[
 L^p(\mathcal M,\tau)\text{ has the UAP}
 \quad\Longleftrightarrow\quad
 L^p(\mathcal M,\tau)\text{ has the column--row property}.
\]
For example, in the hyperfinite \(\mathrm{II}_1\) factor \(R\), the
trace-preserving conditional expectations onto \(M_2^{\otimes n}\)
are finite-rank contractions converging pointwise in \(L^p(R)\)
to the identity \cite[pp.~1465--1466, 1490--1491]{px}. Hence
\(L^p(R)\) has the metric approximation property. Since \(R\) is a
type \(\mathrm{II}_1\) factor, Theorem~\ref{thm:classification} gives
failure of the UAP for \(1\le p<\infty\), \(p\ne2\).

\subsection{Countable discrete groups}\label{sec:groups}

\begin{corollary}\label{cor:groups}
Let \(\Gamma\) be a countable discrete group and let \(1\le p<\infty\),
\(p\ne2\). Then \(L^p(L\Gamma)\), formed with the canonical trace, has the UAP if
and only if \(\Gamma\) contains an abelian subgroup of finite index.
\end{corollary}

\begin{proof}
If \(L^p(L\Gamma)\) has the UAP, Theorem~\ref{thm:classification} makes
\(L\Gamma\) subhomogeneous and hence of type~I. The regular
representation of \(\Gamma\) is therefore of type~I. By the
Thoma--Kaniuth theorem \cite[Thm.~E, (iii)$\Rightarrow$(ii)]{bekka},
\(\Gamma\) has an abelian subgroup of finite index.

Conversely, let \(H\le\Gamma\) be abelian with \([\Gamma:H]=d<\infty\).
By \cite[Thm.~E, (ii)$\Rightarrow$(i), (iv)]{bekka}, \(\Gamma\) is
of type~I and each irreducible unitary representation \(\pi\) is
finite dimensional. Choose a common eigenvector \(v\ne0\) for \(\pi(H)\)
and left coset representatives \(g_1,\ldots,g_d\). By irreducibility,
\[
 \mathcal H_\pi=\operatorname{span}\{\pi(g)v:g\in\Gamma\}
 =\operatorname{span}\{\pi(g_j)v:1\le j\le d\},\qquad
 \dim\mathcal H_\pi\le d .
\]
In the type~I Plancherel decomposition of \(L\Gamma\), the factor
belonging to \(\pi\) is \(B(\mathcal H_\pi)\)
\cite[Thm.~A and Sec.~1]{bekka}. Thus all matrix degrees are at
most \(d\). Theorem~\ref{thm:classification} gives the UAP.
\end{proof}

In particular, \(L^p(L\mathbb F_2)\) fails the UAP for
\(1\le p<\infty\), \(p\ne2\). Indeed, every finite-index subgroup of
\(\mathbb F_2=\langle a,b\rangle\) contains some positive powers of
\(a\) and \(b\), and those powers do not commute.

\subsection{Failure of the ordinary approximation property}

\begin{corollary}\label{cor:ap}
For every \(1<p<\infty\) with \(p\ne2\), there exists a finite von
Neumann algebra \(N\) with separable predual and a faithful normal
tracial state \(\tau\) such that \(L^p(N,\tau)\) fails the
Banach-space approximation property.
\end{corollary}

\begin{proof}
Fix \(1<p<\infty\), \(p\ne2\). Suppose that \(L^p(N,\tau)\) has the AP
for every finite \(N\) with separable predual. Under this assumption,
Claims I--III in the proof of \cite[Thm.~2.19, pp.~284--285]{jnrx}
give the metric approximation property for \(L^p(M)\) for every von
Neumann algebra \(M\). By Raynaud's ultrapower theorem \cite[Thm.~3.1]{raynaud},
each Banach ultrapower of \(S_p\) is isometric to such an \(L^p(M)\).
Thus all ultrapowers of \(S_p\) would have the metric approximation
property, as in Claim~IV of \cite[p.~285]{jnrx}. Heinrich's criterion
\cite[Thm.~9.1(ii), p.~97]{heinrich} would then give the UAP for
\(S_p\), contrary to Theorem~\ref{thm:main}.
\end{proof}

The exponent range in Corollary~\ref{cor:ap} is optimal within
\(1<p<\infty\). For every finite von Neumann algebra \(N\),
\(L^2(N)\) is a Hilbert space and its finite-dimensional orthogonal
projections give the metric approximation property.

\begin{corollary}\label{cor:rigid-factor}
For every countable \(D\subset(1,\infty)\setminus\{2\}\), there is a
\(\mathrm{II}_1\) factor \(M\) with separable predual and property~\((T)\)
such that
\begin{equation}\label{eq:rigid-factor}
 \operatorname{Out}(M)=\{e\},\qquad \mathcal F(M)=\{1\},\qquad
 L^p(M)\text{ fails the ordinary AP for every }p\in D .
\end{equation}
Here \(\operatorname{Out}(M)\) is the outer automorphism group and
\(\mathcal F(M)\) is the fundamental group.
\end{corollary}

\begin{proof}
It suffices to consider nonempty \(D\). For each \(p\in D\), take
\((N_p,\tau_p)\) from Corollary~\ref{cor:ap}. Choose \(w_p>0\)
with \(\sum_{p\in D}w_p=1\). The von Neumann algebra
\[
 N=\prod_{p\in D}N_p,\qquad
 \tau((x_p)_{p\in D})=\sum_{p\in D}w_p\tau_p(x_p)
\]
has a faithful normal tracial state and separable predual
\(\bigl(\bigoplus_{p\in D}(N_p)_*\bigr)_{\ell_1}\).
By Chifan, Drimbe and Ioana \cite[Thm.~A(2)]{cdi}, \((N,\tau)\)
embeds trace-preservingly into a \(\mathrm{II}_1\) factor \(M\) with
separable predual, property~\((T)\) and the two groups in
\eqref{eq:rigid-factor}. Identify \(N\) with its image and let
\(\mathbb E_N:M\to N\) be the trace-preserving conditional expectation.
It extends to a contraction on each \(L^p\)
\cite[Sec.~7, pp.~1490--1491]{px}.

Fix \(p\in D\) and let \(\iota_p:N_p\to N\) be the inclusion into
the \(p\)-coordinate. For \(x\in L^p(N_p,\tau_p)\) and \(y\in L^p(M)\), put
\begin{equation}\label{eq:factor-retraction}
 J_p x=w_p^{-1/p}\iota_p(x),\qquad
 Q_p y=w_p^{1/p}(\mathbb E_N y)_p .
\end{equation}
The identity \(\|z\|_{L^p(N,\tau)}^p
=\sum_{q\in D}w_q\|z_q\|_{L^p(N_q,\tau_q)}^p\) shows that
\(J_p\) is isometric and \(Q_p\) is contractive. For every compact
\(C\subset L^p(N_p,\tau_p)\) and every finite-rank operator \(T\) on
\(L^p(M)\), \eqref{eq:factor-retraction} gives
\[
 \sup_{x\in C}\|Q_pTJ_px-x\|_p
 \le\sup_{y\in J_p(C)}\|Ty-y\|_p,\qquad
 \operatorname{rank}(Q_pTJ_p)\le\operatorname{rank}T .
\]
The AP on \(L^p(M)\) would therefore imply the AP on \(L^p(N_p,\tau_p)\),
contrary to Corollary~\ref{cor:ap}.
\end{proof}

In particular, taking \(D=\{2+1/n:n\ge1\}\) in
Corollary~\ref{cor:rigid-factor} gives a single \(M\) for which
\(L^{2+1/n}(M)\) fails the ordinary AP for every \(n\ge1\).

\subsection{Schatten ultrapowers and two questions of Junge, Nielsen, Ruan and Xu}

Fix \(1<p<\infty\). In Problems~2.20(i) and (iv), Junge, Nielsen,
Ruan and Xu ask whether every \(\mathcal{OL}_p\)-space has the CBAP
and whether an ultrapower of \(S_p\) is an \(\mathcal{OL}_p\)-space,
respectively \cite[p.~285]{jnrx}. Here \(\mathcal{OL}_p\) denotes their
class of operator spaces in which every finite-dimensional subspace
is contained in a finite-dimensional subspace at uniformly bounded
cb-distance from some \(L^p(A)\) with \(A\) a finite-dimensional
\(C^*\)-algebra \cite[Lem.~2.1(i)]{jnrx}. They observed that
an affirmative answer to (iv) forces a negative answer to (i) for
\(p>80\) \cite[p.~285]{jnrx}.

\begin{corollary}\label{cor:olp}
For each \(1<p<\infty\), \(p\ne2\), some operator-space ultrapower of
\(S_p\) fails the ordinary bounded approximation property. Consequently,
at each such \(p\), Problems~2.20(i) and (iv) of \cite{jnrx} cannot
both have affirmative answers.
\end{corollary}

\begin{proof}
By Theorem~\ref{thm:main} and Heinrich's ultrapower criterion
\cite[Thm.~9.1(ii), pp.~97--98]{heinrich}, there is a free ultrafilter
\(\mathcal U\) on \(\mathbb N\) for which the Banach ultrapower
\(X=(S_p)_{\mathcal U}\) fails the bounded approximation property.
With its canonical operator-space structure, \(X\) also fails the
CBAP, since \(\|T\|\le\|T\|_{\mathrm{cb}}\) for every completely
bounded map \(T:X\to X\).
If \(X\) is an \(\mathcal{OL}_p\)-space, it is a counterexample to (i).
If \(X\) is not an \(\mathcal{OL}_p\)-space, it is a counterexample to (iv).
\end{proof}

Failure of the CBAP for some Schatten ultrapower also follows from
the higher-rank lattice results. For each \(1<p<\infty\), \(p\ne2\),
choose \(r\) so that \(L^p(L\Gamma)\) fails the CBAP for
\(\Gamma=\mathrm{SL}(r,\mathbb Z)\) \cite[Thm.~1.1]{dlds}.
Reduction modulo primes separates the elements of \(\Gamma\), so
\(\Gamma\) is residually finite and \(L\Gamma\) has the quotient weak
expectation property \cite[Sec.~5]{jr}. Hence \(L^p(L\Gamma)\) is
completely contractively complemented in a Schatten ultrapower
\cite[proof of Thm.~2.17]{jnrx}. A completely bounded approximation
on that ultrapower would restrict, after projection, to one on
\(L^p(L\Gamma)\). Corollary~\ref{cor:olp} obtains the stronger
obstruction of failure of the ordinary bounded approximation property.

\section{An independent column--row norm formula}\label{sec:independent-row}

For \(n\ge1\), equip \(S_p^n\) with the canonical operator space
structure used in \cite[Sec.~2]{rs}. The column--row property of an
operator space means that the norms of the maps from block columns to
block rows are bounded uniformly in their length. Ray and Sarkar
conjectured in 2023 that their lower bound for the following map is
sharp \cite[Rem.~2.7]{rs}. Its norm can be computed by complex
interpolation alone.

\begin{proposition}\label{prop:column-row}
For \(1\le p\le\infty\), let \(T_n:M_{n,1}(S_p^n)\to M_{1,n}(S_p^n)\)
send the block column with entries \(x_1,\ldots,x_n\) to the block
row with the same entries. Its ordinary operator norm is
\begin{equation}\label{eq:column-row}
 \|T_n\|=n^{|1/2-1/p|}.
\end{equation}
\end{proposition}

\begin{proof}
Fix \(n\) and write \(a_p=\|T_n\|\). Ray and Sarkar
\cite[Lem.~2.6 and Rem.~2.8]{rs} give
\[
 a_p\ge n^{|1/2-1/p|}\quad(p\ne2),\qquad a_1,a_\infty\le\sqrt n .
\]
At \(p=2\), let \(x\) be the block column with entries \(x_i\).
The operator Hilbert space formula \cite[(2.2) and Lem.~2.3]{rs} gives
\[
 \|T_nx\|_{M_{1,n}(S_2^n)}^4
 =\sum_{i,j=1}^n|\operatorname{Tr}(x_i x_j^*)|^2
 =\|x\|_{M_{n,1}(S_2^n)}^4 .
\]
Hence \(a_2=1\). Let \(0<\theta<1\) and
\(1/p=(1-\theta)/p_0+\theta/p_1\). For \((r,s)=(n,1)\) and \((1,n)\), Pisier's
rectangular interpolation identity \cite[Sec.~2.7, (2.7.5)]{pisier-os}
and reiteration for the canonical Schatten scale give, isometrically,
\begin{equation}\label{eq:corner-interpolation}
 [M_{r,s}(S_{p_0}^n),M_{r,s}(S_{p_1}^n)]_\theta=M_{r,s}(S_p^n).
\end{equation}
Interpolating the same linear map \(T_n\) on the two couples in
\eqref{eq:corner-interpolation} gives
\[
 a_p\le
 \begin{cases}
 a_1^{\,2/p-1}a_2^{\,2-2/p}\le n^{1/p-1/2}, &1<p<2,\\[2pt]
 a_2^{\,2/p}a_\infty^{\,1-2/p}\le n^{1/2-1/p}, &2<p<\infty .
 \end{cases}
\]
Together with the endpoint estimates and the lower bound, this proves
\eqref{eq:column-row}.
\end{proof}

\end{document}